\documentclass[12pt,twoside]{amsart}
\usepackage[ansinew]{inputenc}

\usepackage{amssymb}
\usepackage{amsmath}
\usepackage{bbm}
\usepackage{mathrsfs}

{\theoremstyle{definition}
\newtheorem{ttt-s}{}[subsection]
\newtheorem{ttt}{}[section]

\newtheorem{defi-s}[ttt-s]{Definition}
\newtheorem{rem-s}[ttt-s]{Remark}
\newtheorem{pr-s}[ttt-s]{Problem}
\newtheorem{al-s}[ttt-s]{Algorithm}

\newtheorem{rem}[ttt]{Remark}
\newtheorem{rems}[ttt]{Remarks}
}

\newtheorem{theo-s}[ttt-s]{Theorem}
\newtheorem{prop-s}[ttt-s]{Proposition}
\newtheorem{flem-s}[ttt-s]{Fundamental Lemma}

\newtheorem{theo}[ttt]{Theorem}
\newtheorem{lem}[ttt]{Lemma}

\newcounter{abc}
\newenvironment{abc}{\begin{list}{\rm \alph{abc}) }{\usecounter{abc} \leftmargin=0.0pt \labelsep=0.0pt \listparindent=0.0pt \labelwidth=0.0pt \parsep=\smallskipamount \itemsep=0.0pt \topsep=0.0pt \partopsep=\smallskipamount}}{\end{list}}

\newcounter{iii}
\newenvironment{iii}{\begin{list}{\rm \roman{iii}) }{\usecounter{iii} \leftmargin=0.0pt \labelsep=0.0pt \listparindent=0.0pt \labelwidth=0.0pt \parsep=\smallskipamount \itemsep=0.0pt \topsep=0.0pt \partopsep=\smallskipamount}}{\end{list}}

\newcounter{III}
\newenvironment{III}{\begin{list}{\rm \Roman{III}. }{\usecounter{III} \leftmargin=0.0pt \labelsep=0.0pt \listparindent=0.0pt \labelwidth=0.0pt \parsep=\smallskipamount \itemsep=0.0pt \topsep=0.0pt \partopsep=\smallskipamount}}{\end{list}}

\usepackage{fancybox}
\newenvironment{quell}{\VerbatimEnvironment
\begin{Tiny}\begin{Verbatim}}{\end{Verbatim}\end{Tiny}}

\newcommand{\bbF}{{\mathbbm F}}
\newcommand{\bbN}{{\mathbbm N}}
\newcommand{\bbQ}{{\mathbbm Q}}
\newcommand{\bbR}{{\mathbbm R}}
\newcommand{\bbZ}{{\mathbbm Z}}

\newcommand{\br}{\discretionary{}{}{}\penalty0\relax}

\newcommand{\lcm}{\mathop{\rm lcm}\limits}
\renewcommand{\gcd}{\mathop{\rm gcd}\limits}

\newcommand{\REDC}{\mathop{\rm REDC}}

\newcommand{\notd}{\mathord{\nmid}}

\def\pmod#1{\nobreak\ifinner\mkern8mu\else\mkern18mu\fi (\text{\rmfamily\upshape mod}\,\,#1)}

\def\rightend#1#2{{%
 \leavevmode\nobreak\hskip .5em plus 1fil
 \penalty600 \hskip 0pt plus -1filll
 \vadjust{}\nobreak\hskip 0pt plus 1filll%
 #1\parfillskip=#2\relax \par}}

\def\eop{\ifmmode\rule[-22pt]{0pt}{1pt}\ifinner\tag*{$\square$}\else\eqno{\square}\fi\else\rightend{$\square$}{0pt}\fi}

\def\rmtext#1{{\text{\rmfamily\mdseries\upshape #1}}}
\def\ittext#1{{\text{\rmfamily\mdseries\itshape #1}}}

\renewcommand{\mod}{\rmtext{ mod }}

\makeatletter
\def\hsmash{\relax 
  \ifmmode\def\next{\mathpalette\mathhsm@sh}\else\let\next\makehsm@sh
  \fi\next}
\def\makehsm@sh#1{\setbox\z@\hbox{#1}\finhsm@sh}
\def\mathhsm@sh#1#2{\setbox\z@\hbox{$\m@th#1{#2}$}\finhsm@sh}
\def\finhsm@sh{\wd\z@\z@ \box\z@}
\makeatother

\renewcommand{\thefootnote}{\fnsymbol{footnote}}
\author[Andreas-Stephan Elsenhans and J\"org Jahnel]{Andreas-Stephan Elsenhans and J\"org Jahnel}
\date{}

\title[The Fibonacci sequence modulo $p^2$]{The Fibonacci sequence modulo $p^2$ -- \\ An investigation by computer for $p < 10^{14}$}

\begin{document}

\begin{abstract}
We show that for primes
$p < 10^{14}$
the period length
$\kappa (p^2)$
of the Fibonacci sequence modulo
$p^2$
is never equal to its period length modulo
$p$.\linebreak[3]
The~investigation involves an extensive search by computer. As an application, we establish the general formula
$\kappa (p^n) = \kappa (p) \cdot p^{n-1}$
for all primes less than
$10^{14}$.
\end{abstract}

\maketitle

\section{Introduction}
\footnotetext[0]{While this work was done the first author was supported in part by a Doctoral Fellowship of the Deutsche Forschungsgemeinschaft~(DFG).}
\footnotetext[0]{The computer part of this work was executed on the Linux~PCs of the Gau\ss\ Laboratory for Scientific Computing at the G\"ottingen Mathematical Institute. Both authors are grateful to Prof.~Y.~Tschinkel for the permission to use these machines as well as to the system administrators for their~support.}

\thispagestyle{empty}

\begin{ttt}
The Fibonacci sequence
$\{ F_k \}_{k \geq 0}$
is defined recursively by
$F_0 = 0$,
$F_1 = 1$,
and
$F_k = F_{k-1} + F_{k-2}$
for
$k \geq 2$.
Modulo some integer
$l \geq 2$,
it must ultimately become periodic as there are only
$l^2$
different pairs of residues
\mbox{modulo~$l$.}
There~is no pre-period since the recursion may be reversed to
$F_{k-2} = F_k - F_{k-1}$.
The~minimal period
$\kappa (l)$
of the Fibonacci sequence
\mbox{modulo~$l$}
is often called the {\em Wall number\/} as its main properties were discovered by~D.~D.~Wall~\cite{Wall}.

Wall's results may be summarized by the theorem below. It shows, in particular, that
$\kappa (l)$
is in general a lot smaller
than~$l^2$.
In fact, one always has
$\kappa (l) \leq 6l$
whereas~equality holds if and only if
$l = 2 \cdot 5^n$
for some
$n \geq 1$.
\end{ttt}

\begin{theo}[{\rm Wall}{}]
\label{Wall}
\begin{abc}
\item
If\/
$\gcd (l_1, l_2) = 1$
then
$\kappa (l_1 l_2) = \lcm (\kappa (l_1), \kappa (l_2))$.

In particular, if\/
$l = \prod_{i=1}^N p_i^{n_i}$
where the
$p_i$
are pairwise different prime numbers
then~$\kappa (l) = \lcm (\kappa (p_1^{n_1}), \ldots, \kappa (p_N^{n_N}))$.

It is therefore sufficient to understand\/
$\kappa$
on prime powers.
\item
$\kappa (2) = 3$
and
$\kappa (5) = 20$.
Otherwise,

$\bullet$
if
$p$
is a prime such that
$p \equiv \pm1 \pmod 5$
then~$\kappa (p) | (p-1)$.

$\bullet$
If
$p$
is a prime such that
$p \equiv \pm2 \pmod 5$
then
$\kappa (p) | (2p+2)$
but~$\kappa (p) \notd (p+1)$.
\item
If
$l \geq 3$
then
$\kappa (l)$
is even.
\item
If
$p$
is prime,
$e \geq 1$,
and
$p^e | F_{\kappa (p)}$
but
$p^{e+1} \notd F_{\kappa (p)}$
then
\begin{equation}
\kappa (p^n) = \left\{
\begin{array}{ll}
\kappa(p)                   & \ittext{for } n \leq e, \\
\kappa(p) \!\cdot\! p^{n-e} & \ittext{for } n > e.
\end{array}
\right.
\end{equation}
\end{abc}
\end{theo}

\section{The Open Problems}
\subsection{The Period Length Modulo a Prime}

\begin{ttt-s}
It is quite surprising that the Fibonacci sequence still keeps secrets. But~there are at least two of them.
\end{ttt-s}

\begin{pr-s}
The first open problem is ``What is the exact value
of~$\kappa (p)$?''.
Equivalently, one should understand precisely the behaviour of the
quotient~$Q$
given by
$\smash{Q (p) := \frac{p-1}{\kappa (p)}}$
for
$p \equiv \pm1 \pmod 5$
and
$\smash{Q (p) := \frac{2(p+1)}{\kappa (p)}}$
for
$p \equiv \pm2 \pmod 5$.
One~might hope for a formula expressing
$Q (p)$
in terms
of~$p$
but, may be, that is too optimistic.
\end{pr-s}

\begin{ttt-s}
It is known that
$Q$
is unbounded. This is an elementary result due to D.~Jarden~\cite[Theorem~3]{Jarden}.

On the other hand,
$Q$
does not at all tend to infinity. If fact, in his unpublished Ph.D.~thesis~\cite{Gottsch}, G.~G\"ottsch computes a~certain average value
of~$\frac1Q$.
To be more precise, under the assumption of the Generalized Riemann Hypothesis, he~proves
$$\sum_{\substack{p \equiv \pm1 \pmod 5 \\ p \leq x, \, p \rmtext{ prime}}} \frac1{Q(p)} = C_1 \frac{x}{\log x} + O \Big( \frac{x \log \log x}{\log^2 x} \Big)$$
where
$\smash{C_1 = \frac{342}{595} \prod_{p \rmtext{ prime}} (1 - \frac{p}{p^3-1}) \approx 0.331\,055\,98}$.
The proof shows as well that the density~of
$\{ p \rmtext{ prime} \mid Q (p) = 1, p \equiv \pm1 \pmod 5 \}$
within the set of all primes is equal
to~$\smash{C_2 = \frac{27}{38} \prod_{p \rmtext{ prime}} (1 - \frac1{p(p-1)}) \approx 0.265\,705\,45}$.

Not assuming any hypothesis, it is still possible to verify that the right hand side constitutes an upper bound. For that, the error term needs to be weakened
to~$\smash{O (\frac{x \log\log\log x}{\log x \log\log x})}$.

For the case
$p \equiv \pm2 \pmod 5$,
G.~G\"ottsch's results are less~strong. Under the assumption of the Generalized Riemann Hypothesis, he establishes the~estimate
$$\sum_{\substack{p \equiv \pm2 \pmod 5 \\ p \equiv 3 \pmod 4 \\ p \leq x, \, p \rmtext{ prime}}} \frac1{Q(p)} \leq C_2 \frac{x}{\log x} + O \Big( \frac{x \log \log \log x}{\log x \cdot \log\log x} \Big)$$
where
$\smash{C_3 = \frac14 \prod_{p \rmtext{ prime}, p \neq 2,5} (1 - \frac{p}{p^3-1}) \approx 0.210\,055\,99}$.
The density of the~set
$\{ p \rmtext{ prime} \mid Q (p) = 1, p \equiv \pm2 \pmod 5, p \equiv 3 \pmod 4 \}$
within the set of all primes is at most
$\smash{C_4 = \frac14 \prod_{p \rmtext{ prime}, p \neq 2,5} (1 - \frac1{p(p-1)}) \approx 0.196\,818\,85}$.
\end{ttt-s}

\begin{ttt-s}
It seems, however, that the inequalities could well be equalities. In~addition, the restriction to primes satisfying
$p \equiv 3 \pmod 4$
might be~irrelevant.

In~fact, we performed a~count for small~primes
$p < 2 \cdot 10^7$
by~computer. Up~to that bound, there are
$317\,687$
prime~numbers such that
$p \equiv \pm2 \pmod 5$
and
$p \equiv 3 \pmod 4$.
At~them, we find
$Q (p) = 1$
exactly
$250\,246$
times which is a relative frequency of
$0.787\,712\,434\;\ldots\; = 4 \cdot 0.196\,928\,108\;\ldots~$.

On the other hand, there are 
$317\,747$
primes
$p$
satisfying
$p \equiv \pm2 \pmod 5$
and
$p \equiv 1 \pmod 4$.
Among them,
$Q (p) = 1$
occurs
$250\,353$
times which is basically the same~frequency as in the
case~$p \equiv 3 \pmod 4$.
\end{ttt-s}

\subsection{The Period Length Modulo a Prime Power}

\begin{pr-s}
\label{prob}
There is another open problem. In fact, there is one question which was left open in the formulation of Theorem~\ref{Wall}: What is the exact value
of~$e$
in dependence
of~$p$?
Experiments for
small~$p$
show
that~$e=1$.
Is this always the~case? In~other words, does one always have
\begin{equation}
\label{phi}
\kappa (p^n) = \kappa (p) \cdot p^{n-1}
\end{equation}
similarly to the famous formula for Euler's
$\varphi$~function?

This is the most perplexing point in D.~D.~Wall's whole study of the Fibonacci sequence modulo~$m$.
For
$p < 10^4$,
it was investigated by help of an electronic computer by Wall in 1960, already.
\end{pr-s}

\begin{ttt-s}
We continued Wall's investigation concerning Problem~\ref{prob} on today's machines. Our main result is Theorem~\ref{main}~below. The~purpose of the present article is to give a description of our approach, particularly of the various algorithms developed and optimizations~used.
\end{ttt-s}

\begin{defi-s}
We call a prime number
$p$
{\em exceptional\/} if equation~(\ref{phi}) is wrong for
some~$n \geq 2$.
\end{defi-s}

\begin{flem-s}
\label{2}
Assume
$p \neq 2$
and
$l$
to be a multiple of
$\kappa (p)$.
Then,~$p^e | F_l$
is sufficient for
$l$
being a period of\/
$\{ F_k \mod p^e \}_{k \geq 0}$,
i.e.~for
$\kappa (p^e) | l$.

{\bf Proof.}
{\em
The claim is that, in our situation,
$F_{l+1} \equiv 1 \pmod {p^e}$
is automatic.

For that, we note that there is the standard~formula
$F_{l+1} F_{l-1} - F_l^2 = (-1)^l = 1$
which we explain in~(\ref{stan}) below. Here, by assumption,
$F_l \equiv 0 \pmod {p^e}$
and, by virtue of the recursion,
$F_{l-1} \equiv F_{l+1} \pmod {p^e}$.
Therefore,
$F_{l+1}^2 \equiv 1 \pmod {p^e}$.

On the other hand, the condition
$\kappa (p) | l$
implies
that~$F_{l+1} \equiv 1 \pmod p$.

As
$p \neq 2$,
Hensel's lemma says that the lift is unique. This 
shows~$F_{l+1} \equiv 1 \pmod {p^e}$
which is our claim.\eop
}
\end{flem-s}\pagebreak[2]

\begin{prop-s}
\label{exc}
Let\/
$p$
be a prime number. Then, the following assertions are~equivalent.
\begin{iii}
\item
$p$
is exceptional,
\item
$F_{\kappa (p)}$
is divisible
by~$p^2$.
\end{iii}
{\bf Proof.}
{\em
``\,i) $\Longrightarrow$ ii)''
Assume, to the contrary, that
$p^2 \notd F_{\kappa (p)}$.
By definition of
$\kappa (p)$,
we know for sure that nevertheless
$p | F_{\kappa (p)}$.
Together, these statements mean, Theorem~\ref{Wall}.d) may be applied for
$e=1$
showing
$\kappa (p^n) = \kappa (p) \cdot p^{n-1}$
for
every~$n \in \bbN$.
This contradicts~i).

``\,ii) $\Longrightarrow$ i)''
We choose the maximal
$e \in \bbN$
such that
$p^e | F_{\kappa (p)}$.
By~assumption,
$e \geq 2$.
Then, Theorem~\ref{Wall}.d) implies
$\kappa (p^2) = \kappa (p)$
which shows equation~(\ref{phi}) to be wrong
for~$n = 2$.
$p$
is exceptional.\eop
}
\end{prop-s}

\begin{prop-s}
\label{ex}
Let\/
$p \neq 2,5$
be a prime number.
\begin{III}
\item
If\/
$p \equiv \pm1 \pmod 5$
then the following assertions are~equivalent.
\begin{iii}
\item
$p$
is exceptional,
\item
$F_{p-1}$
is divisible
by~$p^2$,
\item
For\/
$\smash{r = \frac{1 + \sqrt{5}}2 \in \bbZ / p^2 \bbZ}$
one has\/
$r^{p-1} = 1$.
\end{iii}
\item
If\/
$p \equiv \pm2 \pmod 5$
then the following assertions are equivalent.
\begin{iii}
\item
$p$
is exceptional,
\item
$F_{2p+2}$
is divisible
by~$p^2$,
\item
$\smash{F_{p+1}}$
is divisible
by~$p^2$.
\item
In\/
$\smash{R_p := \bbZ / p^2 \bbZ\, [r] / (r^2 - r - 1)}$
one has\/
$\smash{r^{p+1} = -1}$.
\end{iii}
\end{III}
{\bf Proof.}
{\em
I.
``\,i) $\Longrightarrow$ iii)''
We put
$s = \frac{1 - \sqrt{5}}2$
and use formula~(\ref{Binet})~below. By Proposition~\ref{exc},
$F_{\kappa (p)}$
is divisible
by~$p^2$.
Therefore,
$r^{\kappa (p)} = s^{\kappa (p)} = (r^{\kappa (p)})^{-1} \in (\bbZ / p^2 \bbZ)^*$,
i.e.~$(r^{\kappa (p)})^2 = 1$.
Since~$\kappa (p) | (p-1)$,
we may conclude
$(r^{p-1})^2 = 1$
from~this. On~the other hand, we~know 
$r^{p-1} \equiv 1 \pmod p$
by Fermat's~Theorem. Uniqueness~of Hensel's~lift
implies~$r^{p-1} = 1$.

``\,iii) $\Longrightarrow$ ii)''
We have
$r^{p-1} s^{p-1} = (-1)^{p-1} = 1$.
Thus,
$r^{p-1} = 1$
implies
$s^{p-1} = 1$.
Consequently,
$\smash{(F_{p-1} \mod p^2) = \frac{r^{p-1} - s^{p-1}}{\sqrt{5}}} = 0$
and
$F_{p-1}$
is divisible
by~$p^2$.

``\,ii) $\Longrightarrow$ i)''
As
$(p-1)$
is a multiple of
$\kappa (p)$,
Lemma~\ref{2} may be applied. It~shows
$\kappa (p^2) | (p-1)$.
This contradicts equation~(\ref{phi}) for
$n=2$.
$p$~is~exceptional.

II.
``\,i) $\Longrightarrow$ ii)''
By Proposition~\ref{exc},
$F_{\kappa (p)}$
is divisible
by~$p^2$.
In that situation, Lemma~\ref{2} implies that
$\kappa (p)$
is actually a period of
$\{ F_k \mod p^2 \}_{k \geq 0}$.
By consequence,
$(2p + 2)$
is a period of
$\{ F_k \mod p^2 \}_{k \geq 0}$,
too. This~shows
$p^2 | F_{2p+2}$.

``\,ii) $\Longrightarrow$ iii)''
Since 
$F_{2p+2} = F_{p+1} V_{p+1}$,
all we need is
$p \,\notd\, V_{p+1}$.
This, however, is clear as
$\smash{V_{p+1} = r^{p+1} + \overline{r}^{p+1}} \equiv -2 \pmod p$.

``\,iii) $\Longrightarrow$ iv)''
The assumption implies
$r^{p+1} = \overline{r}^{p+1}$,
i.e.~
$r^{p+1} \in \bbZ / p^2 \bbZ$.
As~$r\overline{r} = -1$,
we may conclude
$(r^{p+1})^2 = r^{p+1} \overline{r}^{p+1} = (-1)^{p+1} = 1$
from~this. Hensel's~lemma implies
$r^{p+1} = -1$
since
$r^{p+1} \equiv -1 \pmod p$
is~known.

``\,iv) $\Longrightarrow$ i)''
$r^{p+1} = -1$
makes sure that
$\smash{(F_{2p+2} \mod p^2) = \frac{(r^{p+1})^2 - (\overline{r}^{p+1})^2}{\sqrt{5}} = 0}$.
As~$(2p + 2)$
is a multiple of
$\kappa (p)$,
Lemma~\ref{2} may be applied. It shows
$\kappa (p^2) | (2p + 2)$.
This contradicts equation~(\ref{phi}) for
$n=2$.
$p$~is~exceptional.
\eop
}
\end{prop-s}

\begin{rem-s}
By Proposition~\ref{ex}, the problem of finding exceptional primes is in perfect analogy to the problem of finding {\em Wieferich primes}.

In~the Wieferich case, one knows
$2^{p-1} \equiv 1 \pmod p$
and would like to understand the set of all primes for which even
$2^{p-1} \equiv 1 \pmod {p^2}$
is~valid. Here,~we know
$F_{\kappa (p)} \equiv 0 \pmod p$
and look for the primes which fulfill
$F_{\kappa (p)} \equiv 0 \pmod {p^2}$.

At least in the case
$p \equiv \pm1 \pmod 5$,
there is, in fact, more than just an analogy. We~consider a particular case of the generalized Wieferich problem
where~$2$
is replaced
by~$r$.
\end{rem-s}

\begin{rem-s}
One might want to put the concept of an exceptional prime into the wider context of algebraic number~theory. We work in the number~field
$\bbQ \bigl( \! \sqrt{5} \bigr)$
in which
$\smash{r = \frac{1+\sqrt{5}}2}$
is a~fundamental~unit.

By analogy, we~could say that an odd~prime
number~$p$
is {\em exceptional for the real quadratic number
field\/}~$\smash{K = \bbQ \bigl( \! \sqrt{d} \bigr)}$~if,
for
$\varepsilon$
a~fundamental unit
in~$K$,
$\varepsilon^{p - 1} \equiv 1 \pmod {p^2}$
when
$\smash{\bigl( \frac{d}{p} \bigr) = 1}$,
$\varepsilon^{2p + 2} \equiv 1 \pmod {p^2}$
when
$\smash{\bigl( \frac{d}{p} \bigr) = -1}$,
or
$\varepsilon^{p (p-1)} \equiv 1 \pmod {p^2}$
in~the ramified
case~$d | p$.
A~congruence
modulo~$p^2$
is, of course, supposed to mean equality in
${\mathscr O}_K / (p^2)$
which, as
$p \neq 2$,
is isomorphic to
$\bbZ [\sqrt{d}] / (p^2) = \bbZ / p^2 \bbZ \, [X] / (X^2 - d)$.
Note~that there is no ambiguity coming from the choice
of~$\varepsilon$
since all exponents are~even.

For many real quadratic number fields it does not require sophisticated programming to find a few exceptional~primes. Below,~we~give the complete list of all exceptional primes
$p < 10^9$
for the fields
$\smash{\bbQ \bigl( \! \sqrt{d} \bigr)}$
where
$d$
is square-free and
up~to~$101$.
Thereby,~primes put in parentheses are those such that
$\bbQ \bigl( \! \sqrt{d} \bigr)$
is ramified
at~$p$.%

{\tt\Tiny
\tabcolsep4.4pt
\begin{center}
\begin{tabular}{|r|l||r|l||r|l|}
\hline
 d & exceptional primes p & d & exceptional primes p & d & exceptional primes p \\
\hline
 2 & 13, 31, 1\,546\,463                   & 17 &                                             & 34 & 37, 547, 4\,733                   \\
 3 & 103                                   & 19 & 79, 1\,271\,731, 13\,599\,893, 31\,352\,389 & 35 & 23, 577, 1\,325\,663              \\
 5 &                                       & 21 & 46\,179\,311                                & 37 & 7, 89, 257, 631                   \\
 6 & (3), 7, 523                           & 22 & 43, 73, 409, 28\,477                        & 38 & 5                                 \\
 7 &                                       & 23 & 7, 733                                      & 39 & 5, 7, 37, 163\,409, 795\,490\,667 \\
10 & 191, 643, 134\,339, 25\,233\,137      & 26 & 2\,683, 3\,967, 18\,587                     & 41 & 29, 53, 7\,211                    \\
11 &                                       & 29 & 3, 11                                       & 42 & (3), 5, 43, 71                    \\
13 & 241                                   & 30 &                                             & 43 & 3, 479                            \\
14 & 6\,707\,879, 93\,140\,353             & 31 & 157, 261\,687\,119                          & 46 & (23)                              \\
15 & (3), 181, 1\,039, 2\,917, 2\,401\,457 & 33 & (3), 29, 37, 6\,713\,797                    & 47 & 5\,762\,437                       \\
\hline\hline
 d & exceptional primes p & d & exceptional primes p & d & exceptional primes p \\
\hline
51 & (3), 5, 37, 4\,831                    & 67 & 3, 11, 953, 57\,301                         &  83 & 3, 19\,699, 2\,417\,377          \\
53 & 5                                     & 69 & (3), 5, 17, 52\,469\,057                    &  85 & 3, 204\,520\,559                 \\
55 & 571                                   & 70 & (5), 59, 20\,411                            &  86 & 1\,231, 5\,779                   \\
57 & 59, 28\,927, 1\,726\,079, 7\,480\,159 & 71 & 67, 2\,953, 8\,863, 522\,647\,821           &  87 & (3), 17, 757, 1\,123             \\
58 & 3, 23, 4\,639, 172\,721, 16\,557\,419 & 73 & 5, 7, 41, 3\,947, 6\,079                    &  89 & 5, 7, 13, 59                     \\
59 & 1\,559, 17\,385\,737                  & 74 & 3, 7, 1\,171                                &  91 & (13), 1\,218\,691                \\
61 &                                       & 77 & 3, 418\,270\,987                            &  93 & (3), 13                          \\
62 & 3, 5, 263, 388\,897                   & 78 & (3), 19, 62\,591                            &  94 & 73                               \\
65 & 1\,327, 8\,831, 569\,831              & 79 & 3, 113, 4\,049, 6\,199                      &  95 & 6\,257, 10\,937                  \\
66 & 21\,023, 106\,107\,779                & 82 & 3, 5, 11, 769, 3\,256\,531, 624\,451\,181   &  97 & 17, 3\,331                       \\
   &                                       &    &                                             & 101 & 7, 19\,301                       \\
\hline
\end{tabular}
\end{center}
}

Among these
$158$
exceptional primes, there are exactly nine for which even the stronger congruence
modulo~$p^3$
is~true. These~are
$p = 3$
for
$d = 29$,
$42$,
$67$,
and~$74$,
$p = 5$
for
$d = 62$,
$73$,
and~$89$,
$p = 17$
for
$d = 69$,
and
$p = 29$
for
$d = 41$.
We~do not observe a congruence of the type above
modulo~$p^4$.
\end{rem-s}\pagebreak[3]

\section{Background}

\begin{ttt}
Part a) of Wall's theorem is trivial.

For the {\bf proof} of b), Binet's formula
\begin{equation}
\label{Binet}
F_k = \frac{r^k - s^k}{\sqrt{5}},
\end{equation}
where
$\smash{r = \frac{1+\sqrt{5}}{2}}$
and
$\smash{s = \frac{1-\sqrt{5}}{2}}$,
is of fundamental importance. It is easily established by induction. If
$p \equiv \pm1 \pmod 5$
then
$5$
is a quadratic residue modulo
$p$
and, therefore,
$\smash{\frac{1 \pm \sqrt{5}}{2} \in \bbF_p}$.
Fermat states their order is a divisor of
$p-1$.

Otherwise, 
$\smash{\frac{1 \pm \sqrt{5}}{2} \in \bbF_{p^2}}$
are elements of
norm~$(-1)$.
As the norm map
$\smash{N \colon \bbF_{p^2}^* \to \bbF_p^*}$
is surjective, its kernel is a group of order 
$\smash{\frac{p^2-1}{p-1}} = p+1$
and
$\# N^{-1} (\{\, 1, -1 \,\}) = 2p+2$.

As
$\bbF_{p^2}^*$
is cyclic, we see that
$N^{-1} (\{\, 1, -1 \,\})$
is even a cyclic group of
order~$2p+2$.
$N (r) = N (s) = -1$
implies that both
$r$
and
$s$
are not contained in its subgroup of index~two. Therefore,
\begin{equation}
r^{p+1} \equiv s^{p+1} \equiv -1 \pmod p.
\end{equation}
From this, we find
$\smash{F_{p+2} \equiv \frac{r^{p+2} - s^{p+2}}{\sqrt{5}} \equiv \frac{-r + s}{\sqrt{5}} \equiv -F_1 \equiv -1 \pmod p}$
which shows
$p+1$
is not a period of
$\{ F_k \}_{k \geq 0}$
modulo~$p$.

c)
In the case
$p \equiv \pm2 \pmod 5$
this follows from~b). It is, however, true in general.

Indeed, for every
$k \in \bbN$,
one has
\begin{eqnarray}
\label{stan}
F_{k+1} F_{k-1} - F_k^2 & = & \nonumber \frac{r^{2k} + s^{2k} - r^{k+1}s^{k-1} - r^{k-1}s^{k+1}}{5} - \frac{r^{2k} + s^{2k} - 2r^ks^k}{5} \\
 & = & \frac{-(-1)^{k-1} (r^2 + s^2) + 2 (-1)^k}{5} \\
 & = & \nonumber (-1)^k
\end{eqnarray}
as
$rs = -1$
and
$r^2 + s^2 = 3$.
On the other hand,
$$F_{\kappa (l)+1} F_{\kappa (l)-1} - F_{\kappa (l)}^2 \equiv 1 \cdot 1 - 0^2 \equiv 1 \pmod l.$$
As
$l \geq 3$
this implies
$\kappa (l)$
is even.

For d), it is best to establish the following
$p$-uplication
formula first.
\end{ttt}

\begin{lem}[{\rm Wall}{}]
\label{1}
One has
\begin{equation}
F_{pk} = \frac{1}{2^{p-1}} \sum_{\substack{{j=1} \\ j \, {\rm odd}}}^p \binom{p}{j} 5^\frac{j-1}{2} F_k^j V_k^{p-j}.
\end{equation}
Here,
$\{ V_k \}_{k \geq 0}$
is the\/ {\em Lucas sequence} given by
$V_0 = 2$,
$V_1 = 1$,
and
$V_k = V_{k-1} + V_{k-2}$
for
$k \geq 2$.

{\bf Proof.}
{\em
Induction shows
$V_k = r^k + s^k$.
Having that in mind, it is easy to calculate as follows.
$$F_{pk} = \frac{(r^k)^p - (s^k)^p}{\sqrt{5}} = \frac{(\frac{V_k + \sqrt{5} F_k}{2})^p - (\frac{V_k - \sqrt{5} F_k}{2})^p}{\sqrt{5}}.$$
The assertion follows from the Binomial Theorem.\eop
}
\end{lem}\pagebreak[1]

\begin{ttt}
The fundamental Lemma~\ref{2} allows us to prove d) for
$p \neq 2$
in a somewhat simpler manner than D.~D.~Wall did it in~\cite{Wall}.

First, we note that for
$n \leq e$,
\ref{2}~implies
$\kappa (p^n) | \kappa (p)$.
However, divisibility the other way round is obvious.

For
$n \geq e$,
by Lemma~\ref{2},
it is sufficient to prove
$\nu_p (F_{\kappa (p) \cdot p^{n-e}}) = n$,
i.e.~that
$p^n | F_{\kappa (p) \cdot p^{n-e}}$
but
$p^{n+1} \notd F_{\kappa (p) \cdot p^{n-e}}$.
Indeed, the first divisibility implies
$\kappa (p^n) | \kappa (p) \cdot p^{n-e}$
while the second, applied for
$n-1$
instead of
$n$,
yields
$\kappa (p^n) \notd \kappa (p) \cdot p^{n-e-1}$.
The result follows as
$\kappa (p) | \kappa (p^n)$.

For
$\nu_p (F_{\kappa (p) \cdot p^{n-e}}) = n$,
we proceed by induction, the case
$n=e$
being known by assumption. One has
$$F_{\kappa (p) \cdot p^{n-e+1}} = \textstyle{\frac{1}{2^{p-1}} p F_{\kappa (p) \cdot p^{n-e}} V_{\kappa (p) \cdot p^{n-e}}^{p-1} + \frac{1}{2^{p-1}} \sum\limits_{\substack{{j=3} \\ j \, {\rm odd}}}^p \binom{p}{j} 5^\frac{j-1}{2} F_{\kappa (p) \cdot p^{n-e}}^j V_{\kappa (p) \cdot p^{n-e}}^{p-j}}.$$
In the second term, every summand is divisible by 
$\smash{F_{\kappa (p) \cdot p^{n-e}}^3}$,
i.e.~by~$p^{3n}$.
The claim would follow if we knew
$p \,\notd V_{\kappa (p) \cdot p^{n-e}}$.
This, however, is easy as there is the formula
\begin{equation}
V_l = F_{l-1} + F_{l+1}
\end{equation}
which implies
$V_l \equiv 2 \pmod p$
for~$l$~any
multiple
of~$\kappa (p)$.
\end{ttt}

\begin{ttt}
For
$p=2$,
as always, things are a bit more complicated. We still have
$\kappa (2^n) = 3 \cdot 2^{n-1}$.
However, for
$n \geq 2$,
one has
$2^{n+1} | F_{3 \cdot 2^{n-1}}$
for which there is no analogue in the
$p \neq 2$
case. On the other hand,
$\nu_2 (F_{3 \cdot 2^{n-1} + 1} - 1) = n$
which is sufficient for our assertion.

The duplication formula provided by Lemma~\ref{1} is
\begin{equation}
F_{2k} = F_k V_k = F_k (F_{k-1} + F_{k+1}) = F_k^2 + 2 F_k F_{k-1}.
\end{equation}
As
$F_6 = 8$,
a~repeated application of this formula shows
$2^{n+1} | F_{3 \cdot 2^{n-1}}$
for
every~$n \geq 2$.

We further claim
$F_{2k+1} = F_k^2 + F_{k+1}^2$.
Indeed, this is true for
$k=0$
as
$1 = 0^2 + 1^2$
and we proceed by induction as follows:
\begin{equation}
\begin{split}
F_{2k+3} = F_{2k+1} + F_{2k+2} = F_k^2 + F_{k+1}^2 + F_{k+1}^2 + 2 F&_{k+1} F_k = \\
= F_{k+1}^2 + (&F_k + F_{k+1})^2 = F_{k+1}^2 + F_{k+2}^2.
\end{split}
\end{equation}
The assertion 
$\nu_2 (F_{3 \cdot 2^{n-1} + 1} - 1) = n$
is now easily established by induction. We~note that
$F_7 = 13 \equiv 1 \pmod 4$
but the same is no~longer true
modulo~$8$.
Furthermore,
$F_{3 \cdot 2^n + 1} = F_{3 \cdot 2^{n-1}}^2 + F_{3 \cdot 2^{n-1} + 1}^2$
where the first summand is even divisible by
$2^{2n+2}$.
The~second one is congruent
to~$1$
modulo
$2^{n+1}$,
but not modulo
$2^{n+2}$,
by consequence of the induction hypothesis.
\end{ttt}

\section{A heuristic argument}

\begin{ttt}
We expect that there are infinitely many exceptional
primes~for~$\bbQ \bigl( \! \sqrt{5} \bigr)$.

Our reasoning for this is as follows.
$p | F_{\kappa (p)}$
is known by definition of
$\kappa (p)$.
Thus, for any individual
prime~$p$,
$(F_{\kappa (p)} \mod p^2)$
is one residue out
of~$p$~possibilities.
If~we were allowed to assume equidistribution then we could conclude that
$p^2 | F_{\kappa (p)}$
should occur with a~``probability''
of~$\smash{\frac1p}$.
Further, by~\cite[Theorem~5]{RS},
$$\log \log N + A - \frac1{2\log^2 N} \leq \sum_{\substack{p \rmtext{ prime} \\ p \leq N}} \frac1p \leq \log \log N + A + \frac1{2\log^2 N},$$
at least for
$N \geq 286$.
Here,~$A \in \bbR$
is Mertens'~constant which is given by
$$A = \gamma + \sum_{p \rmtext{ prime}} \left[ \frac1p + \log \Big( 1 - \frac1p \Big) \right] = 0.261\,497\,212\,847\,642\,783\,755\;\ldots$$
whereas
$\gamma$
denotes the Euler-Mascheroni constant.

This means that one should expect around
$\log \log N + A$
exceptional primes less
than~$N$.
\end{ttt}

\begin{ttt}
On the other hand,
$p^3 | F_{\kappa (p)}$
should occur only a few times or even not at all. Indeed, if we assume equidistribution again, then for any individual
prime~$p$,
$p^3 | F_{\kappa (p)}$
should happen with a~``probability''
of~$\frac1{p^2}$.
However,
$$\sum_{\substack{p=2 \\ p \rmtext{ prime}}}^\infty \frac1{p^2} = 0.452\,247\,420\,041\,065\,498\,506\;\ldots\;~.$$
is a convergent series.
\end{ttt}

\begin{rem}
It is, may be, of interest that, for any
exponent~$n \geq 2$,
one has the equality
$\sum_{p \rmtext{ prime}} \frac1{p^n} = \sum_{k=1}^\infty \frac{\mu(k)}k \log \zeta (nk)$
where the right hand converges a lot faster and may be used for~evaluation. This equation results from the Moebius inversion formula and Euler's formula
$\log \zeta (nk) = \;\sum\limits_{\hsmash{p \rmtext{ prime}}} - \log (1 - \frac1{p^{nk}}) = \sum_{j=1}^\infty \frac1j \;\sum\limits_{\hsmash{p \rmtext{ prime}}}\; \frac1{p^{jnk}}$.
\end{rem}

\begin{ttt}
\label{main}
We carried out an extensive search for exceptional primes but, unfortunately, we had no success and our result is negative.\medskip

{\bf Theorem.}
{\em
There are no exceptional
primes~$p < 10^{14}$.

Down the earth, this means that one has
$\kappa (p^n) = \kappa (p) \cdot p^{n-1}$
for
every~$n \in \bbN$
and all
primes~$p < 10^{14}$.
}
\end{ttt}

\section{Algorithms}
\setcounter{ttt-s}{0}

\begin{ttt-s}
We worked with two principally different types of algorithms. First, in the
$p \equiv \pm1 \pmod 5$
case, it is possible to compute
$(r^{p-1} \mod {p^2})$.
A second and more complete approach is to compute
$(F_{p-1} \mod {p^2})$
in the
$p \equiv \pm1 \pmod 5$
case and
$(F_{2p+2} \mod {p^2})$
or
$(F_{p+1} \mod {p^2})$
in the case
$p \equiv \pm2 \pmod 5$.
\end{ttt-s}

\begin{rem-s}
In the case
$p \equiv \pm2 \pmod 5$,
$p \neq 2$,
exceptionality is equivalent to
$r^{2p+2} \equiv 1 \pmod {p^2}$.
Unfortunately, an approach based on that observation turns out to be impractical as it involves the calculation of a modular power in
$R_p = \bbZ / p^2 \bbZ \bigl[ \sqrt{5} \bigr] = \bbZ \bigl[ \sqrt{5} \bigr] / (p^2)$
in a situation where
$\sqrt{5} \not\in \bbZ / p^2 \bbZ$.
In comparison with
$\bbZ / p^2 \bbZ$,
multiplication in
$R_p$
is a lot slower, at least in our (naive) implementations. This puts a modular powering operation
in~$R_p$
out of competition with a direct approach to compute 
$F_{2p+2} \rmtext{ (or } F_{p+1})$
modulo~${p^2}$.
\end{rem-s}

\subsection{Algorithms based on the computation of $\sqrt{5}$}

\begin{ttt-s}
\label{alg}
If
$p \equiv \pm1 \pmod 5$
then one may routinely compute
$(r^{p-1} \mod {p^2})$.
The algorithm should consist of four steps.
\begin{iii}
\item
Compute the square root of
$5$
in
$\bbZ / p \bbZ$.
\item
Take the Hensel's lift of this root to
$\bbZ / p^2 \bbZ$.
\item
Calculate the golden ratio
$\smash{r := \frac{1+\sqrt{5}}{2}} \in \bbZ / p^2 \bbZ$.
\item
Use a modular powering operation to find
$(r^{p-1} \mod {p^2})$.
\end{iii}

We call algorithms which follow this strategy {\em algorithms powering the golden~ratio.}

Here, the final steps iii) and iv) are not critical at all. For~iii), it is obvious that this is a simple calculation while for~iv), carefully optimized modular powering operations are available. Further, ii)~can be effectively done as
$r^2 \equiv 5 \pmod p$
implies
$\smash{w := r - \frac{r^2 - 5}{p} \cdot (\frac{1}{2r} \mod p) \cdot p}$
is a square root
of~$5$
modulo~$p^2$.
Thus, the most expensive operation should be a run of Euclid's extended algorithm in~order to
find~$(\frac{1}{2r} \mod p)$.

In fact, there is a way to avoid even~this. We~first
calculate~$\smash{\frac15 \in \bbF_p}$.
This is easier than an arbitrary division in residues
modulo~$p$.
We may put
$\frac15 := \frac{4p+1}5$
if~$p \equiv 1 \pmod 5$
and
$\frac15 := \frac{p+1}5$
if~$p \equiv -1 \pmod 5$.
Then,~a~representative~$v$
of~$(\frac1r \mod p)$
can be computed as
$v = r \cdot \frac15$.
We get away with one integer~division and one~multiplication.
\end{ttt-s}

\begin{ttt-s}
Thus, the most interesting point is~i), the computation
of~$\sqrt{5} \in \bbF_p$.
In~general, there is a beautiful algorithm to find square roots modulo a prime~number due to Shanks~\cite[Algorithm~1.5.1]{Cohen}. We~implemented this algorithm but let it finally run only in the
$p \equiv 1 \pmod 8$
case. If
$p \not\equiv 1 \pmod 8$
then there are direct formulae to compute the square root
of~$5$
which turn out to work faster.

If~$p \equiv 3 \pmod 4$
then one may simply put
$\smash{w := (5^{\frac{p+1}4} \mod p)}$
to find a square root
of~$5$
by one modular powering~operation.

If
$p \equiv 5 \pmod 8$
then one may put
\begin{equation}
\label{F1}
w := (5^{\frac{p+3}8} \mod p)
\end{equation}
as long as
$5^{\frac{p-1}4} \equiv 1 \pmod p$
and
\begin{equation}
\label{F2}
w := (10 \cdot 20^{\frac{p-5}8} \mod p)
\end{equation}
if
$5^{\frac{p-1}4} \equiv -1 \pmod p$.
Note that
$5$
is a quadratic residue
modulo~$p$.
Hence, we always have 
$\smash{5^{\frac{p-1}4} \equiv \pm1 \pmod p}$.

For sure,
$\smash{(5^{\frac{p-1}4} \mod p)}$
can be computed using a modular powering operation. In~fact, we implemented an algorithm doing that and let it run through the
interval~$[10^{12}, 5 \cdot 10^{12}]$.

However,
$(5^{\frac{p-1}4} \mod p)$
is nothing but a quartic residue symbol. For that reason, there is an actually faster algorithm which we obtained by an approach using the law of biquadratic~reciprocity.
\end{ttt-s}

\begin{prop-s}
\label{viert}
Let
$p$
be a prime number such that
$p \equiv 5 \pmod 8$
and
$p \equiv \pm1 \pmod 5$
and let
$p = a^2 + b^2$
be its (essentially unique) decomposition into a~sum of two squares.
\begin{abc}
\item
Then,
$a$
and
$b$
may be normalized such that
$a \equiv 3 \pmod 4$
and
$b$
is even.
\item
Assume\/
$a$
and\/
$b$
are normalized as described in a). Then, there are only the following eight possibilities.
\begin{iii}
\item
$a \equiv 3, 7, 11, \ittext{or } 19 \pmod {20}$
and\/
$\,b \equiv 10 \pmod {20}$.
\newline
In this case,
$\smash{5^{\frac{p-1}4} \equiv 1 \pmod p}$,
i.e.~$5$~is~a
quartic residue
modulo~$p$.
\item
$a \equiv 15 \pmod {20}$
and\/
$b \equiv 2, 6, 14, \ittext{or } 18 \pmod {20}$.
\newline
Here,
$\smash{5^{\frac{p-1}4} \equiv -1 \pmod p}$,
i.e.~$5$~is~a
quadratic but not a quartic residue
modulo~$p$.
\end{iii}
\end{abc}

{\bf Proof.}
{\em
a)
As~$p$
is odd, among the integers
$a$
and
$b$
there must be an even and an odd one. We choose
$b$
to be even and force
$a \equiv 3 \pmod 4$
by replacing
$a$
by~$(-a)$,
if~necessary.

b)
We first observe that
$a^2 \equiv 1 \pmod 8$
forces
$b^2 \equiv 4 \pmod 8$
and
$b \equiv 2 \pmod 4$.
Then, we realize that one of the two numbers
$a$
and
$b$
must be divisible
by~$5$.
Indeed, otherwise we had
$a^2, b^2 \equiv \pm1 \pmod 5$
which does not allow
$a^2 + b^2 \equiv \pm1 \pmod 5$.
Clearly,
$a$
and
$b$
cannot be both divisible
by~$5$.

If
$a$
is divisible
by~$5$
then 
$a \equiv 3 \pmod 4$
implies
$a \equiv 15 \pmod {20}$.
$b \equiv 2 \pmod 4$
and
$b$~not
divisible
by~$5$
yield the four possibilities stated. On the other hand, if
$b$~is
divisible
by~$5$
then 
$b \equiv 2 \pmod 4$
implies
$b \equiv 10 \pmod {20}$.
$a \equiv 3 \pmod 4$
and
$a$~not
divisible
by~$5$
show there are precisely the four possibilities listed.

For the remaining assertions, we first note that
$(5^{\frac{p-1}4} \mod p)$
tests whether
$x^4 \equiv 5 \pmod p$
has a solution
$x \in \bbZ$,
i.e.~whether
$5$
is a quartic residue
modulo~$p$.
By~\cite[Lemma~9.10.1]{IR}, we know
$$(5^{\frac{p-1}4} \mod p) = \chi_{a + bi} (5)$$
where
$\chi$
denotes the quartic residue symbol. The law of biquadratic~reciprocity\linebreak[1] \cite[Theorem~9.2]{IR} asserts
$$\chi_{a + bi} (5) = \chi_5 (a + bi).$$
For that, we note explicitly that
$a + bi \equiv 3 + 2i \pmod 4$,
$5 \equiv 1 \pmod 4$,
and
$\smash{\frac{N (5) - 1}{4} = 6}$
is even. Let us now
compute~$\chi_5 (a + bi)$:
\begin{align*}
\chi_5 (a + bi) & = \chi_{-1 + 2i} (a + bi) \cdot \chi_{-1 - 2i} (a + bi) \notag \\
                & = \overline{\chi_{-1 - 2i} (a - bi)} \cdot \chi_{-1 - 2i} (a + bi) \displaybreak[0] \notag \\
                & = \Big( \frac{a + \frac{b}2}{5} \Big) \cdot \chi_{-1 - 2i} (a - bi) \cdot \chi_{-1 - 2i} (a + bi) \\
                & = \Big( \frac{a + \frac{b}2}{5} \Big) \cdot \chi_{-1 - 2i} (p). \notag
\end{align*}
Here, the first equation is the definition of the quartic residue symbol for composite elements while the second is~\cite[Proposition~9.8.3.c)]{IR}.

For the third equation, we observe that
$\chi_{-1 - 2i} (a - bi)$
is either
$\pm1$
or
$\pm{i}$.
By~simply omitting the complex conjugation, we would make a sign error if and only if
$\chi_{-1 - 2i} (a - bi) = \pm{i}$.
By~\cite[Lemma~9.10.1]{IR}, this means exactly that
$a - bi$
defines, under the identification
$2i = -1$,
not even a quadratic residue
modulo~$5$.\linebreak[2]
Therefore, the correction factor is
$\smash{(\frac{a + \frac{b}2}{5})}$.
The final equation follows from~\cite[Proposition~9.8.3.b)]{IR}.

We note that, by virtue of~\cite[Lemma~9.10.1]{IR},
$\chi_{-1 - 2i} (p)$
tests
whether~$p$
is a quartic residue
modulo~$5$
or not.
As~$p$
is for sure a quadratic residue, we may write
\begin{eqnarray*}
\chi_{-1 - 2i} (p) = \left\{
\begin{array}{ll}
\phantom{-}1 & \rmtext{ if } p \equiv \phantom{-}1 \pmod 5, \\
-1           & \rmtext{ if } p \equiv -1 \pmod 5
\end{array}
\right.
\end{eqnarray*}
or, if we want,
$\chi_{-1 - 2i} (p) = (p \mod 5)$.

The eight possibilities could now be inspected one after the other. A more conceptual argument works as follows. In case~i), we have
$$\Big( \frac{a + \frac{b}2}{5} \Big) = \Big( \frac{a}5 \Big) = (a^2 \mod 5) = (a^2 + b^2 \mod 5) = (p \mod 5).$$
Therefore,
$(5^{\frac{p-1}4} \mod p) = 1$.
On the other hand, in case~ii),
\begin{equation*}
\begin{split}
\Big( \frac{a + \frac{b}2}{5} \Big) = \Big( \frac{\frac{b}2}5 \Big) = \Big( \frac{b^2}4 \mod 5 \Big) = (-b^2 \mod 5) = (-a^2 - b^2 &\mod 5) = \\
&= -(p \mod 5).
\end{split}
\end{equation*}
Hence,
$(5^{\frac{p-1}4} \mod p) = -1$.\eop
}
\end{prop-s}

\begin{ttt-s}
Although we are not going to make use of it, let us state the complementary result
for~$p \equiv 1 \pmod 8$.\smallskip

{\bf Proposition.}
{\em
Let
$p$
be a prime such that
$p \equiv 1 \pmod 8$
and
$p \equiv \pm1 \pmod 5$
and let
$p = a^2 + b^2$
be its (essentially unique) decomposition into a~sum of two squares.
\begin{abc}
\item
Then,
$a$
and
$b$
may be normalized such that
$a \equiv 1 \pmod 4$
and
$b$
is even.
\item
Assume\/
$a$
and\/
$b$
are normalized as described in a). Then, there are only the following eight possibilities.
\begin{iii}
\item
$a \equiv 1, 9, 13, \ittext{or } 17 \pmod {20}$
and\/
$\,b \equiv 0 \pmod {20}$.
\newline
In this case,
$\smash{5^{\frac{p-1}4} \equiv 1 \pmod p}$,
i.e.~$5$~is~a
quartic residue
modulo~$p$.
\item
$a \equiv 5 \pmod {20}$
and\/
$b \equiv 4, 8, 12, \ittext{or } 16 \pmod {20}$.
\newline
Here,
$\smash{5^{\frac{p-1}4} \equiv -1 \pmod p}$,
i.e.~$5$~is~a
quadratic but not a quartic residue
modulo~$p$.
\end{iii}
\end{abc}
}

{\bf Proof.}
a)
As~$p$
is odd, among the integers
$a$
and
$b$
there must be an even and an odd one. We choose
$b$
to be even and force
$a \equiv 1 \pmod 4$
by replacing
$a$
by~$(-a)$,
if~necessary.

b)
We first observe that
$a^2 \equiv 1 \pmod 8$
forces
$b^2 \equiv 0 \pmod 8$
and
$4 | b$.
Then, we realize that one of the two numbers
$a$
and
$b$
must be divisible
by~$5$.
Indeed, otherwise we had
$a^2, b^2 \equiv \pm1 \pmod 5$
which does not allow
$a^2 + b^2 \equiv \pm1 \pmod 5$.
Clearly,
$a$
and
$b$
cannot be both divisible
by~$5$.

If
$a$
is divisible
by~$5$
then 
$a \equiv 1 \pmod 4$
implies
$a \equiv 5 \pmod {20}$.
$4 | b$
and
$b$~not
divisible
by~$5$
yield the four possibilities stated. On the other hand, if
$b$~is
divisible
by~$5$
then 
$4 | b$
implies
$b \equiv 0 \pmod {20}$.
$a \equiv 1 \pmod 4$
and
$a$~not
divisible
by~$5$
show there are precisely the four possibilities listed.

The proof of the remaining assertions works exactly in the same way as the proof of Proposition~\ref{viert}~above. We note explicitly that
$a + bi \equiv 1 \pmod 4$
makes sure that the law of biquadratic reciprocity may be~applied.\eop
\end{ttt-s}

\begin{ttt-s}
As the transformation
$a \mapsto -a$
does not affect any of the three statements below, we may formulate the following theorem. Actually, this~is the result we need for the~application.\smallskip\pagebreak[4]

{\bf Theorem.}
{\em
Let
$p$
be a prime number such that
$p \equiv 1 \pmod 4$
and
$p \equiv \pm1 \pmod 5$
and let
$p = a^2 + b^2$
be its decomposition into a~sum of two squares. We normalize
$a$
and
$b$
such that
$a$
is~odd and
$b$
is~even. Then, the following three statements are~equivalent.
\begin{iii}
\item
$5$
is a quartic residue
modulo~$p$.
\item
$b$
is divisible
by~$5$.
\item
$a$
is not divisible
by~$5$.
\end{iii}
}
\end{ttt-s}

\begin{rem-s}
We note that the restrictions
on~$p$
exclude only trivial~cases.
If~$p \not\equiv \pm1 \pmod 5$
then~$5$
is not even a quadratic residue
modulo~$p$.
If~$p \equiv 3 \pmod 4$
then every quadratic residue is automatically a quartic~residue.
\end{rem-s}

\begin{al-s}
The {\em square sum sieve algorithm\/} for prime numbers~$p$
such that
$p \equiv 21, 29 \pmod {40}$
runs as follows.

We investigate a rectangle
$[N_1, N_2] \times [M_1, M_2]$
of numbers. We will go through the rectangle row-by-row in the same way as the electron beam goes through a screen.
\begin{abc}
\item
We add
$0$,
$1$,
$2$,
or
$3$
to
$M_1$
to make sure
$M_1 \equiv 2 \pmod 4$.
Then, we let
$b$
go from
$M_1$
to
$M_2$
in steps of length four.
\item
For a fixed
$b$
we sieve the odd numbers in the interval
$[N_1, N_2]$.

Except for the odd case that
$l | a, b$
which we decided to ignore as the density of these pairs is not too high,
$l | a^2 + b^2$
implies that
$(-1)$
is a quadratic residue
modulo~$l$,
i.e.~we need to sieve only by the
primes~$l \equiv 1 \pmod 4$.

For each
such~$l$
which is below a certain limit we cross out all those
$a$
such that
$a \equiv \pm v_l b \pmod l$.
Here,
$v_l$
is a square root
of~$(-1)$
modulo~$l$,
i.e.~$v_l^2 \equiv -1 \pmod l$.
For practical application, this requires that the square roots
of~$(-1)$
modulo the relevant primes have to be pre-computed and stored in an array once and for~all.
\item
For the remaining pairs
$(a, b)$,
we compute
$p = a^2 + b^2$
and do steps i) through~iv) from~\ref{alg}. In step~i),
if~$b$
is divisible
by~$5$
then we use formula~(\ref{F1}) to compute the square root
of~$5$
modulo~$p$.
Otherwise, we use formula~(\ref{F2}).
\end{abc}
\end{al-s}

\begin{ttt-s}
In practice, we ran the square sum sieve algorithm on the rectangles
$[0,4\,000\,000] \times [1\,580\,000,4\,000\,000]$
and
$[1\,580\,000,4\,000\,000] \times [0,1\,580\,000]$,
thereby capturing every prime
$p \in [5 \cdot 10^{12}, 1.6 \cdot 10^{13}]$
such that
$p \equiv 21, 29 \pmod {40}$
plus several~others.

In fact, on the second rectangle we ran a modified version, the {\em inverted square sum sieve,} where the two outer loops are reversed. That means, we~let
$a$
go through the odd numbers in
$[N_1, N_2]$
in the very outer loop. This has some advantage in speed as longer intervals are sieved at once. In other words, we go through the rectangle column-by-column.

We implemented the square sum sieve algorithms in C using the mpz~functions of GNU's GMP package for arithmetic on long integers. On a single 1211\,MHz Athlon processor, the computations for the first rectangle took around 22\,days of CPU time. The computations for the smaller second rectangle were finished after nine~days.
\end{ttt-s}

\begin{ttt-s}
For primes
$p$
such that
$p \equiv 3 \pmod 4$
and
$p \equiv \pm1 \pmod 5$,
the~formula
$\smash{w := (5^{\frac{p+1}4} \mod p)}$
for the square root
of~$5$
makes things a lot easier. Instead of the square sum sieve we implemented the sieve of Eratosthenes. Caused by the limitations of main memory in today's PCs, we could actually sieve intervals of only about
$250\,000\,000$~numbers
at once. For each such interval the remainders of its starting~point have to be computed (painfully) by explicit~divisions.
\end{ttt-s}

\begin{al-s}
More precisely, the {\em algorithm powering the golden ratio for primes\/}
$p \equiv 11, 19 \pmod {20}$
runs as follows.

We investigate an interval
$[N_1, N_2]$.
We~assume that
$N_2 - N_1$
is divisible
by~$5 \cdot 10^9$
and that
$N_1$
is divisible
by~$20$.
\begin{abc}
\item
We let an integer~variable
$i$
count from
$0$~to~$\frac{N_2 - N_1}{5 \cdot 10^9} - 1$.
\item
For fixed
$i$
we work on the
interval~$I = [N_1 + 5 \cdot 10^9 \cdot i, N_1 + 5 \cdot 10^9 \cdot (i + 1)]$.
For each
prime~$l$
which is below a certain limit, we compute
$(N_1 + 5 \cdot 10^9 \cdot i \mod l)$.
Then, we cross out
all~$p \in I$,
$p \equiv 11 \; (\rmtext{or } 19) \mod 20$
which are divisible
by~$l$.
\item
For the remaining
$p \in I$,
$p \equiv 11 \; (\rmtext{or } 19) \mod 20$
we do steps i) through iv) from~\ref{alg}. In step~i), we use the
formula~$\smash{w := (5^{\frac{p+1}4} \mod p)}$
to compute the square root
of~$5$
modulo~$p$.
\end{abc}
\end{al-s}

\begin{ttt-s}
In practice, we ran this algorithm in~order to test all prime numbers
$p \in [10^{12}, 4 \cdot 10^{13}]$
such that
$p \equiv 11 \pmod {20}$
or
$p \equiv 19 \pmod {20}$.
It was implemented in C using the mpz~functions of the GMP package.

Later, when testing primes
above~$10^{13}$,
we used the low level mpn~functions for long natural numbers. In~particular, we implemented a modular powering function which is hand-tailored for numbers of the considered size. It~uses the left-right
base~$2^3$
powering algorithm~\cite[Algorithm~1.2.3]{Cohen} and the sliding window improvement from~mpz\_powm.

Having done all these optimizations, work on the test interval
$[4 \cdot 10^{13}, 4 \cdot 10^{13} + 5 \cdot 10^9]$
of~$250\,000\,000$~numbers
$p$
such that
$p \equiv 11 \pmod {20}$,
among them
$19\,955\,067$
primes, lasted 7:50~Minutes CPU~time on a 1211\,MHz Athlon processor.
\end{ttt-s}

\begin{ttt-s}
Similarly, for prime~numbers
$p$
satisfying the simultaneous congruences
$p \equiv 1 \pmod 8$
and
$p \equiv \pm1 \pmod 5$,
we implemented Shanks' algorithm~\cite[Algorithm~1.5.1]{Cohen}
to compute the square root
of~$5$
modulo~$p$.
\end{ttt-s}\pagebreak[4]

\begin{al-s}
More precisely, the {\em algorithm powering the golden ratio for primes\/}
$p \equiv 1, 9 \pmod {40}$
runs as follows.

We investigate an interval
$[N_1, N_2]$.
We~assume that
$N_2 - N_1$
is divisible
by~$10^{10}$
and that
$N_1$
is divisible
by~$40$.
\begin{abc}
\item
We let an integer~variable
$i$
count from
$0$~to~$\frac{N_2 - N_1}{10^{10}} - 1$.
\item
For fixed
$i$
we work on the
interval~$I = [N_1 + 10^{10} \cdot i, N_1 + 10^{10} \cdot (i + 1)]$.
For each
prime~$l$
which is below a certain limit, we compute
$((N_1 + 10^{10} \cdot i) \mod l)$.
Then, we cross out
all~$p \in I$,
$p \equiv 1 \; (\rmtext{or } 9) \pmod {40}$
which are divisible
by~$l$.
\item
For the remaining
$p \in I$,
$p \equiv 1 \; (\rmtext{or } 9) \pmod {40}$
we do steps i) through iv) from~\ref{alg}. In step~i), we use Shanks'~algorithm to compute the square root
of~$5$
modulo~$p$.
\end{abc}
\end{al-s}

\begin{ttt-s}
We~ran this algorithm on the
interval~$[10^{12}, 4 \cdot 10^{13}]$.
After all optimizations, the test interval
$[4 \cdot 10^{13}, 4 \cdot 10^{13} + 10^{10}]$
of
$250\,000\,000$~numbers
$p$
such that
$p \equiv 1 \pmod {40}$,
among them
$19\,954\,152$
primes, could be searched through on a 1211\,MHz Athlon processor in 10:30~Minutes CPU~time.

This is quite a lot more in comparison with the algorithm for
$p \equiv 11 \pmod {20}$
or~$p \equiv 19 \pmod {20}$.
The~difference comes entirely from the more complicated procedure to
compute~$\sqrt{5} \in \bbF_p$.
\end{ttt-s}

\begin{rem-s}
At a certain moment, such a running time was no~longer found~reasonable. A direct computation of the Fibonacci numbers could be done as~well. After several optimizations of the code of the direct~methods, it turned out that only the
$3$~mod~$4$~case
could still compete with them. We discuss the direct~methods in the subsection~below.
\end{rem-s}

\subsection{Algorithms for a direct computation of Fibonacci numbers}

\begin{al-s}
A nice algorithm for the fast computation of a Fibonacci number is presented in O.~Forster's book~\cite{Forster}. It is based on the formulae
\begin{equation}
\begin{split}
F_{2k-1} & = F_k^2 + F_{k-1}^2, \\
F_{2k}   & = F_k^2 + 2F_k F_{k-1}.
\end{split}
\end{equation}
and works in the spirit of the left-right binary powering algorithm using~bits.

Our adaption uses modular~operations
modulo~$p^2$ 
instead of integer operations. An implementation in O.~Forster's Pascal-style multi~precision interpreter language ARIBAS looks like this.\pagebreak[4]

\begin{quell}
(*------------------------------------------------------------------*)
(*
** Schnelle Berechnung der Fibonacci-Zahlen mittels der Formeln
**      fib(2*k-1) = fib(k)**2 + fib(k-1)**2
**      fib(2*k)   = fib(k)**2 + 2*fib(k)*fib(k-1)
**
** Dabei werden alle Berechnungen mod m durchgeführt
*)
function fib(k,m : integer): integer;
var
    b, x, y, xx, temp: integer;
begin
    if k <= 1 then return k end;
    x := 1; y := 0;
    for b := bit_length(k)-2 to 0 by -1 do
        xx := x*x mod m;
        x := (xx + 2*x*y) mod m;
        y := (xx + y*y) mod m;
        if bit_test(k,b) then
            temp := x;
            x := (x + y) mod m;
            y := temp;
        end;
    end;
    return x;
end.
\end{quell}

\begin{quell}
(** ein systematischer Versuch**)
function test() : integer
var
    p,r,r1 : integer;
    ptest  : boolean;
begin
    for p := 90000000001 to 95000000001 by 2 do
        if (p mod 10000) = 1 then
            writeln("getestete Zahl: ", p);
        end;
        ptest := rab_primetest(p);
        if (ptest = true) then
            if ((p mod 5 = 2) or (p mod 5 = 3)) then
                r := fib(2*p+2,p*p);
            else
                r := fib(p-1,p*p);
            end;
            if (r <= 30000000000000000) then
                r1 := r div p;
                writeln(p," ist eine interessante Primzahl. Quotient ", r1);
            end;
        end;
    end;
    return(0);
end.
\end{quell}

A call to {\tt fib(k,m)} computes
$(F_k \mod m)$.
{\tt test} is the main function. {\tt test()} executes an outer loop which contains a Rabin-Miller composedness test. For~a
pseudo~prime~$p$,
it uses the function fib to compute
$(F_{p-1} \mod p^2)$
or
$(F_{2p+2} \mod p^2)$.
As~these are divisible
by~$p$
we output the quotient~instead. Note that in order to limit the output size we actually write an output only when the quotient is rather~small.
\end{al-s}

\begin{ttt-s}
ARIBAS is fast enough to ensure that this algorithm could be run from
$p = 7$
up~to~$10^{11}$.
We~worked on ten~PCs in parallel for five~days. That~was our first bigger computing project concerning this problem. It showed that no exceptional primes
$p < 10^{11}$
do~exist, thereby a establishing a lightweight~version of Theorem~\ref{main}.
\end{ttt-s}

\begin{ttt-s}
The running time made it clear that we had approached to the limits of an interpreter~language. For a systematic test of larger prime numbers, the algorithm was ported to~C. For the arithmetic on long integers we used the mpz~functions of GMP. After~only one further optimization, the integration of a version of the sieve of Eratosthenes, the
interval~$[10^{11}, 10^{12}]$
could be attacked. A test interval of
$250\,000\,000$~numbers
was dealt with on a 1211\,MHz Athlon processor in around 40~Minutes CPU~time. Again, we did parallel computing on ten~PCs. The search through
$[10^{11}, 10^{12}]$
was finished in less than five~days.
\end{ttt-s}

\begin{ttt-s}
For the
interval~$[10^{12}, 10^{13}]$,
the methods which compute
$\sqrt{5} \in \bbF_p$
and square the golden ratio were introduced as they were faster than our implementation of O.~Forster's algorithm at that time. For this reason, only the
case~$p \equiv \pm2 \pmod 5$
was done by Forster's algorithm. It~took us around 20~days on ten~PCs.
\end{ttt-s}\pagebreak[1]

\section{Optimizations}
\subsection{Sieving}

\begin{ttt-s}
Near
$10^{14}$,
one of about~$32$
numbers is~prime. We work in a fixed prime residue~class
modulo~$10$,
$20$,
or~$40$
but still, only one of about 13~numbers is prime. We~feel that the computations of
$(F_{p \pm 1} \mod p^2)$
should take the main part of the running time of our~programs. Our goal is, therefore, to rapidly exclude (most of) the non-primes from the list and then to spend most of the time on the remaining~numbers.

There are various methods to generate the list of all primes within an~interval. Unfortunately, this section of our code is not as harmless as one could hope~for. In~fact, for an individual
number~$p$,
one might have the idea to decide whether it is probably prime by computing
$(F_{p \pm 1} \mod p)$.
That is the Fibonacci composedness~test. It~would, unfortunately, not reduce our computational load a lot as it is almost as complex as the main~computation. This~clearly indicates the problem that the standard ``pseudo primality tests'' which are designed to test individual numbers are not well suited for our purposes. In~this~subsection, we will explain what we did instead in~order to speed up this~part of the~program.
\end{ttt-s}

\begin{ttt-s}
Our first programs in ARIBAS in fact used the internal primality~test to check each number in the interval individually. At the ARIBAS~level, this is optimal because it involves only one instruction for the~interpreter.

When we migrated our programs to C, using the GMP~library, we first tried the~same. We~used the function mpz\_probab\_prime with one repetition for every number to be tested. It~turned out that this program was enormously inefficient. It took about 50~per~cent of the running~time for primality testing and 50~per~cent for the computation of Fibonacci~numbers. However, it could easily be tuned by a naive implementation of the sieve of Eratosthenes in intervals of
length~$1\,000\,000$.

We first combined sieving by small primes and the mpz\_probab\_prime function because sieving by huge primes is slow. This made sure that the computation of Fibonacci~numbers took the major part of the running~time. However, mpz\_probab\_prime is not at all intended to be combined with a sieve. In fact, it checks divisibility by small primes once more. Thus, an optimization of the code for the Fibonacci~numbers reversed the relation again. It became necessary to carry~out a further optimization of the generation of the list of primes. We decided to abandon all pseudo~primality tests. Further, we enlarged the length of the array of up~to
$250\,000\,000$
numbers to minimize the number of~initializations.

In principle, the sieve works as follows. Recall that we used different algorithms for the computation of the Fibonacci~numbers, depending on the residue~class
of~$p$
modulo~$10$,
$20$,
or~$40$.
This leads to a sieve in which the number 
$$S (i) := \rmtext{starting point} + \rmtext{residue} + \rmtext{modulus} \cdot i$$
is represented by array
position~$i$.
Since all our moduli are divisible by 2 and 5 we do no~longer sieve by these two~numbers.

Such a sieve is still easy to use. Given a
prime~$p \neq 2, 5$,
one has to compute the array
index~$i_0$
of the first number which is divisible
by~$p$.
Then, one can cross out the numbers at the indices
$i_0, i_0 + p, i_0 + 2p, \;\ldots\;$
until the end of the sieve is~reached.
\end{ttt-s}

\begin{ttt-s}
{\bf Optimization for the Cache Memory.}
An array of the size above fits into the memory of today's PCs but it does not fit into the~cache. Thus, the speed-limiting part is the transfer between CPU and memory. Sieving by big primes is like a random access to single bytes. The
memory~manager has to transfer one block to the cache~memory, change one byte, and then transfer the whole block back to the memory. This is the limiting bottleneck.

To avoid this problem as far as possible, we built a {\em two stage sieve.}

In the first stage, we sieve by the first
$25\,000$,
the ``small'', primes. For that, we~divide the sieve further into segments of
length~$30\,000$.
These two constants were found to be optimal in practical~tests. They are heavily machine~dependent.

The first stage is now easily explained. In a first step, we sieve the first segment by all small primes. Then, we sieve the second segment by all small primes. We~continue in that way until the end of the sieve is~reached.

In the second stage, we work with all relevant ``big'' primes on the complete sieve, as~usual.\pagebreak[3]

The result of this strategy is a sieve whose segments fit into the machine's~cache. Thus, the speed of the first sieve stage is the speed of the cache, not the speed of the memory. The speed of the second stage is limited by the initialization. 

On our machines the two stage sieve is twice as fast as the ordinary sieve.
\end{ttt-s}

\begin{ttt-s}
The choice of the prime limit for sieving is a point of interest, too. As~we search for one very particular example, it would do no harm if, from to time, we test a composite
number~$p$
for
$p^2 | F_{p\pm1}$.
When the computer would tell us
$p^2$
divides
$F_{p\pm1}$
which, in fact, it never did then it would be easy to do a reliable primality~test.

As long as we sieve by small primes, it is clear that lots of numbers will be crossed out in a short time and this will reduce the running time as it reduces the number of times the actual computation of
$(F_{p\pm1} \mod p^2)$
is called. Afterwards, when we sieve by larger primes, the situation is no~longer that clear. We will often cross out a number repeatedly which was crossed out already before. This means, it can happen that further sieving costs actually more time than it saves.

Our tests show nevertheless that it is best to sieve {\em almost\/} till to the square root of the numbers to be tested. We introduced an automatic choice of the variable prime\_limit as
$\smash{\frac{\sqrt{p}}{\log{\sqrt{p}}}}$
which means we sieve by the first
$\smash{[\frac{\sqrt{p}}{\log{\sqrt{p}}}]}$
primes.
Here,~$p$
means the first prime of the interval we want to go~through.
\end{ttt-s}

\begin{ttt-s}
Another optimization was done by looking at the prime~three. Every~third number is crossed out when sieving by this prime and, when sieving by a bigger prime, every third step hits a number which is divisible by~three and already crossed~out.

Thus, we can work more efficiently as follows.
Let~$p$
be a prime bigger than~three and coprime to the modulus. We
compute~$i_0$,
the first index of a number divisible
by~$p$.
Then, we calculate the remainder of the corresponding number modulo~three.
If~it is~zero then we skip
$i_0$
and continue with
$i_0 := i_0 + p$.  
Now,
$i_0$
corresponds to the first number in the sieve which is divisible
by~$p$
but not by~three. Thus, we~must~cross~out
$i_0, i_0 + p, i_0 + 3 p, i_0 + 4 p, i_0 + 6p, \;\ldots\;$
or
$i_0, i_0 + 2p, i_0 + 3p, i_0 + 5p, i_0 + 6p, \;\ldots\;$
depending~on whether
$i_0 + 2p$
corresponds to a number which is divisible by~three or~not.
\end{ttt-s}

\subsection{The Montgomery Representation}

\begin{ttt-s}
The algorithms for the computation of Fibonacci numbers
modulo~$m$
explained so far spend the lion's share of their running time on the divisions
by~$m$
which occur as the final steps of modular operations such as
\mbox{\tt x := (xx + 2*x*y) mod m}.
Unfortunately, on today's PC processors, divisions are by far slower than multiplications or even~additions.

An ingenious method to avoid most of the divisions in a modular powering operation is due to P.~L.~Montgomery~\cite{Montgomery}. We use an adaption of Montgomery's method to O.~Forster's algorithm which works as~follows.

Let
$R$
be the smallest positive integer which does not fit into one machine~word. That will normally be a power of two. On our machines,
$R = 2^{32}$.
Recall that all operations on unsigned integers in~C are automatically modular operations
modulo~$R$.
We choose some
exponent~$n$
such that the modulus
$m = p^2$
fulfills
$m \leq \frac{R^n}5$.
In~our situation
$p < 10^{14}$,
therefore
$\smash{m = p^2 < 10^{28} < \frac{2^{96}}5}$,
such that
$n=3$
will be~sufficient.

Instead of the variables
$x, y, \;\ldots\, \in \bbZ / m \bbZ$,
we work with their {\em Montgomery representations\/}
$x_M, y_M, \;\ldots\, \in \bbZ$.
These numbers are not entirely unique but bound to be integers from the
interval~$[0, \frac{R^n}5)$
fulfilling
$x_M \equiv R^n x \pmod m$.
This means that modular divisions still have to be done in some initialization step, one for each variable that is initially there, but these turn out to be the only divisions we are going to~execute!

A modular operation, for example
$x := ((x^2 + 2xy) \mod m)$,
is translated into its {\em Montgomery counterpart.} In the example this is
$$x_M := \Big( \frac{x_M^2 + 2 x_M y_M}{R^n} \mod m \Big).$$
We see here that
$\smash{x_M, y_M < \frac{R^n}5}$
implies
$\smash{x_M^2 + 2 x_M y_M < 3 \cdot \frac{R^{2n}}{25}}$.
An inspection of O.~Forster's algorithm shows that we always have to compute
$(\frac{A}{R^n} \mod m)$
for
some~$\smash{A < 5 \cdot \frac{R^{2n}}{25} = \frac{R^{2n}}5}$.
$$A \mapsto \Big( \frac{A}{R^n} \mod m \Big)$$
is Montgomery's
$\REDC$
function. It~occurs everywhere in the algorithm where normally a reduction
modulo~$m$,
i.e.~$A \mapsto (A \mod m)$,
would be~done.

This looks as if we had not won anything. But, in fact, we won a lot as for computer hardware it is much easier to compute
$(\frac{A}{R^n} \mod m)$,
which is a ``reduction from below'',
than~$(A \mod m)$
which is a ``reduction from above'' and usually involves trial~divisions.

Indeed,
$A$~fits
into
$2n$
machine words. It has
$2n$
so-called {\em limbs}. The rightmost, i.e.~the~least
significant,
$n$
of those have to be transformed into zero by adding some suitable multiple
of~$m$.
Then, these
$n$
limbs may simply be~omitted.

Which multiple
of~$m$
is the suitable one that erases the rightmost limb
$A_0$
of~$A$?
Well,~$q \cdot m$
for
$q := (-A_0 \cdot m^{-1} \mod R)$
will do. This operation is in fact an ordinary multiplication of unsigned integers in~C as
$(-A_0)$
on unsigned integers means
$(R - A_0)$
and multiplication is automatically
modulo~$R$.
We add
$q \cdot m$
to
$A$
and remove the last limb. This procedure of transforming the rightmost machine~word
of~$A$
into zero
and removing it has to be repeated
$n$~times.

Still,
$m$
needs to be inverted modulo
$R = 2^{32}$.
The naive approach for this would be to use Euclid's extended algorithm which, unfortunately, involves quite a number of divisions. At least, we observe that it is necessary to do this only once, not
$n$~times
although there are
$n$~iterations.
However, for the purpose of inverting an odd number
modulo~$2^{32}$,
there exists a very elegant and highly efficient C~macro in GMP, named~modlimb\_invert. It uses a table of the modular inverses of all odd integers
modulo~$2^8$
and then executes two Hensel's lifts in a row. Note that, if
$i \cdot n \equiv 1 \pmod N$
then
$(2i - i^2 \cdot n) \cdot n \equiv 1 \pmod {N^2}$.
We observe that, in this particular case, we need no division for the Hensel's~lift.

What is the size of the representative
of~$(\frac{A}{R^n} \mod m)$
found? We~have
$A < \frac{R^{2n}}5$.
We~add to that less~than
$R^n m$
and divide
by~$R^n$.
Thus, the representative is less~than
$$\frac{\frac{R^{2n}}5 + R^n m}{R^n} = \frac{R^n}5 + m.$$
We want
$\REDC (A) < \frac{R^n}5$,
the same inequality we have for all variables in Montgomery representation. To reach that, we may now simply
subtract~$m$
in the case we found an
outcome~$\geq \frac{R^n}5$.
(This is the point where we use
$m \leq \frac{R^n}5$.)

Our version of
$\REDC$
looks as follows. In order to optimize for speed, we designed it as a C~macro, not as a~function.

\begin{quell}
#define REDC(mp, n, Nprim, tp)                                     \
do {                                                               \
 mp_limb_t cy;                                                     \
 mp_limb_t qu;                                                     \
 mp_size_t j;                                                      \
                                                                   \
 for (j = 0; j < n; j++) {                                         \
  qu = tp[0] * Nprim;                                              \
  /* q = tp[0]*invm mod 2^32. Reduktion mod 2^32 von selber! */    \
  cy = mpn_addmul_1 (tp, mp, n, qu);                               \
  mpn_incr_u (tp + n, cy);                                         \
  tp++;                                                            \
 }                                                                 \
                                                                   \
 if (tp[n - 1] >= 0x33333333)                 /* 2^32 / 5. */      \
  mpn_sub_n (tp, tp, mp, n);                                       \
} while(0);
\end{quell}

It is typically invoked as \mbox{{\tt REDC (m, REDC\_BREITE, invm, ?);},} with various variables in the place of the~?, after invm is set by \mbox{{\tt modlimb\_invert (invm, m[0]);}} and \mbox{{\tt invm = -invm;}.} Up to now, we always had REDC\_BREITE~=~3.

At the very end of our algorithm we find
$(F_k)_M$,
the desired Fibonacci number in its Montgomery representation. To convert back, we just need one more call to
$\REDC$.
Indeed,
$$(F_k \mod m) = \bigg( \frac{F_k R^n}{R^n} \mod m \bigg) = \bigg( \frac{(F_k)_M}{R^n} \mod m \bigg) = \REDC ((F_k)_M).$$
Further,
$(F_k)_M < \frac{R^n}5$
implies
$$\REDC ((F_k)_M) < \frac{\frac{R^n}5 + R^n \cdot m}{R^n} = \frac15 + m,$$
i.e.~$\REDC ((F_k)_M) \leq m$.

We note explicitly that there is quite a dangerous trap at this point. The
residue~$0$,
the one we are in fact looking for, will not be reported
as~$0$
but
as~$m$.
We work around this by outputting residues of small {\em absolute\/} value.
If~$(r \mod m)$
is found and
$r$~is
not below a certain output limit then
$m - r$
is computed and compared with that limit.
\end{ttt-s}

\begin{rem-s}
The integration of the Montgomery representation into our algorithm allowed us to avoid practically all the divisions. This caused a stunning reduction of the running time to about one third of its original value.
\end{rem-s}\pagebreak[1]

\subsection{Other Optimizations}

\begin{ttt-s}
We introduced several other optimizations. One, which is worth a mention, is the integration of a pre-computation for the first seven binary digits
of~$p$.
Note, if we let
$p$~go
linearly through a large interval then its first seven digits will change very slowly. This means, as a study of our algorithm for the computation of
$(F_p \mod p^2)$
shows, that the same first seven steps will be done again and again. We avoid this and do these steps once, as a pre-computation. As
$10^{14}$
consists
of~$47$~binary
digits this saves about
$14$~per~cent
of the running time.

Of course,
$p$
is not a constant for the outer loop of our program and its first seven binary digits are only almost constant. One needs to watch out for the moment when the seventh digit
of~$p$~changes.
\end{ttt-s}

\begin{ttt-s}
Another improvement by a few per~cent was obtained through the switch to a different algorithm for the computation of the Fibonacci numbers. Our~hand-tailored approach computes the
$k$-th
Fibonacci~number~$F_k$
simultaneously with the
$k$-th
Lucas
number~$V_k$.
It~is based on the formulae
\begin{equation}
\label{hand}
\begin{split}
F_{2k}   & = F_k V_k, \\                         
V_{2k}   & = V_k^2 + 2 (-1)^{k+1}, \\
F_{2k+1} & = \frac{F_k V_k + V_k^2}2 + (-1)^{k+1}, \\
V_{2k+1} & = F_{2k+1} + 2 F_k V_k.
\end{split}
\end{equation}
This is faster than the algorithm explained above as it involves only one multiplication and one squaring operation instead of one multiplication and two squaring operations. It seems here that the number of multiplications and the number of squaring operations determine the running time. Multiplications~by~two are not counted as multiplications as they are simple bit shifts. Bit~shifts and additions are a lot faster than multiplications while a squaring operation costs about two thirds of what a multiplication costs.

From that point of view there should exist an even better algorithm. One can make use of the formulae
\begin{equation}
\label{gmp}
\begin{split}
F_{2k+1} & = 4 F_k^2 - F_{k-1}^2 + 2 (-1)^k, \\
F_{2k-1} & = F_k^2 + F_{k-1}^2, \\
F_{2k}   & = F_{2k+1} - F_{2k-1}
\end{split}
\end{equation}
which we found in the GMP source code. If we meet a bit which is set then we continue with
$F_{2k+1}$
and
$F_{2k}$.
Otherwise, with
$F_{2k}$
and
$F_{2k-1}$.

Here, there are only two squaring operations involved and no multiplications, at~all. This should be very hard to beat. Our tests, however, unearthed that the program made from~(\ref{gmp}) ran approximately ten per cent slower than the program made from~(\ref{hand}). For~that reason, we worked finally with~(\ref{hand}). Nevertheless, we expect that for larger
numbers~$p$,
in a situation where additions and bit~shifts contribute even less proportion to the running time, an algorithm using~(\ref{gmp})~should actually run faster. It is possible that this is the case from the moment on that
$p^2 > 2^{96}$
does no~longer fit into three limbs but occupies~four.
\end{ttt-s}

\begin{ttt-s}
Some other optimizations are of a more practical nature. For example, instead of GMP's mpz~functions we used the low level mpn~functions for long natural numbers. Further, we employed some internal GMP~low~level functions although this is not recommended by the GMP documentation.

The point is that the size of the numbers appearing in our calculations is a-priori known to us and basically always the same. When, for example, we multiply two numbers, then it does not make sense always to check whether the base case multiplication, the Karatsuba scheme, or the FFT~algorithm will be fastest. In our case, mpn\_mul\_basecase is always the fastest of the three, therefore we call it~directly.
\end{ttt-s}

\subsection{The Performance Finally Achieved}

\begin{ttt-s}
As a consequence of all the optimizations described, the CPU~time it took our program to test the~interval
$[4 \cdot 10^{13}, 4 \cdot 10^{13} + 2.5 \cdot 10^9]$
of~$250\,000\,000$~numbers
$p$
such that
$p \equiv 3 \pmod {10}$,
among them
$19\,955\,355$
primes, was reduced to~8:08~Minutes. Sieving is done in the first 24~seconds.

The tests were made on a~1211\,MHz Athlon processor. For comparison, on a~1673\,MHz~Athlon processor we test the same interval in around 6:30~Minutes and on a~3\,GHz~Pentium~4 processor in around 5:30~Minutes. (This relatively poor running time might partially be due to the fact that we carried out our trial runs on Athlon~processors.)
\end{ttt-s}

\begin{ttt-s}
{\bf The Main Computational Undertaking.}
In a project of somewhat larger scale, we ran the optimized algorithm on all primes~$p$
in the
interval~$[10^{13}, 10^{14}]$
such that
$p \equiv \pm2 \pmod 5$.
Further, as the methods which start with the computation
of~$\sqrt{5} \in \bbF_p$
are no~longer faster, we ran it, too, on all prime numbers
$p \in [4 \cdot 10^{13}, 10^{14}]$
such that
$p \equiv \pm1 \pmod 5$
and on all primes
$p \in [1.6 \cdot 10^{13}, 4 \cdot 10^{13}]$
such that
$p \equiv 5 \pmod 8$
and
$p \equiv \pm1 \pmod 5$.

Altogether, this means that we fully tested the whole
interval~$[10^{13}, 10^{14}]$.
To~do this took us around 820~days of CPU time. The computational work was done in parallel on up to 14~PCs from July till October~2004.
\end{ttt-s}

\section{Output Data}

\begin{ttt}
{\bf Computer Proof.}
Neither our earlier computations for
$p < 10^{13}$
nor the more recent ones for the
interval~$[10^{13}, 10^{14}]$
detected any exceptional primes. As we covered the intervals systematically and tested each individual prime, this establishes the fact that for all prime
numbers~$p < 10^{14}$
one has
$p^2 \notd F_{\kappa (p)}.$
There are no exceptional primes below that limit. Theorem~\ref{main} is verified.
\end{ttt}

\begin{ttt}
{\bf Statistical Observations.}
We do never find
$(F_{p \pm 1} \mod p^2) = 0$.
Does~that mean, we have found some evidence that our assumption, the residues
$(F_{p \pm 1} \mod p^2)$
should be equidistributed in
$\{\, 0, p, 2p, \;\ldots\;, p^2 - p \,\}$,
is wrong?\linebreak[1] Actually, it does not. Besides the fact that the value zero does not occur, all other reasonable statistical quantities seem to be well within the expected range.

Indeed, a typical piece of our output data looks as follows.

\begin{quell}
Durchsuche Fenster mit Nummer 34304.
Beginne sieben.
Restklassen berechnet.
Beginne sieben mit kleinen Primzahlen.
Sieben mit kleinen Primzahlen fertig.
Fertig mit sieben.
Initialisiere
x mit 110560307156090817237632754212345,
y mit 247220362414275519277277821571239
und vorz mit 1.
 10786 	 Quotient 1912354 mit p := 85760594147971.
 10787 	 Quotient 1072750 mit p := 85760627258851.
 10788 	 Quotient -1617348 mit p := 85760847493241.
 10789 	 Quotient -3142103 mit p := 85761104075891.
Initialisiere
x mit 178890334785183168257455287891792,
y mit 400010949097364802732720796316482
und vorz mit -1.
 10790 	 Quotient -9341211 mit p := 85761921174961.
Fertig mit Fenster mit Nummer 34304.
Durchsuche Fenster mit Nummer 34305.
Beginne sieben.
Restklassen berechnet.
Beginne sieben mit kleinen Primzahlen.
Sieben mit kleinen Primzahlen fertig.
Fertig mit sieben.
Initialisiere
x mit 178890334785183168257455287891792,
y mit 400010949097364802732720796316482
und vorz mit -1.
 10791 	 Quotient 3971074 mit p := 85763512710481.
 10792 	 Quotient 2441663 mit p := 85764391244491.
Fertig mit Fenster mit Nummer 34305.
\end{quell}

To make the output easier to understand we do not print
$(F_{p \pm 1} \mod p^2)$
which is automatically divisible
by~$p$
but
$R (p) := (F_{p \pm 1} \mod p^2) / p \in \bbZ / p \bbZ$.
Such~a quotient may be as large
as~$p$.
We output only those which fall
into~$(-10^7, 10^7)$
which is very small in comparison
to~$p$.

The data above were generated by a process which had started at
$8 \cdot 10^{13}$
and worked on the primes
$p \equiv 1 \pmod {10}$.
Till~$8.5765 \cdot 10^{13}$,
it~found
$10\,792$
primes~$p$
such that
$R (p) = (F_{p-1} \mod p^2) / p \in (-10^7, 10^7)$.

On the other hand, assuming equidistribution we would have predicted to find such a particularly small quotient for~around
\begin{eqnarray*}
\smash{\sum_{\substack{p = 8 \cdot 10^{13} \\ p \equiv 1 \pmod {10} \\ p \rmtext{ prime}}}^{\hsmash{8.5765 \cdot 10^{13}}}} \!\! \frac{2 \!\cdot\! 10^7 - 1}p & \approx & (2 \!\cdot\! 10^7 - 1) \!\cdot\! \frac1{\varphi (10)} \!\cdot\! (\log (\log (8.5765 \!\cdot\! 10^{13})) - \log (\log (8 \!\cdot\! 10^{13}))) \\
 & =       & \frac{2 \!\cdot\! 10^7 - 1}4 \cdot (\log (\log (8.5765 \cdot 10^{13})) - \log (\log (8 \cdot 10^{13}))) \\
 & \approx & 10\,856.330
\end{eqnarray*}
primes which is astonishingly close to the reality.

Among the
$10\,792$
small quotients found within this interval, the absolutely smallest one is
$R (82\,789\,107\,950\,701) = -42$.
We~find
$1\,074$~quotients
of absolute value less than~$1\,000\,000$,
$98$~quotients
of absolute value less than~$100\,000$,
and
$10$
of absolute value less than~$10\,000$.
These are, besides the one above,
\begin{align*}
R (80\,114\,543\,961\,461) & = -2437,                            \\
R (80\,607\,583\,847\,341) & = -6949,                            \\
R (80\,870\,523\,194\,401) & = -5751,           \displaybreak[0] \\
R (81\,232\,564\,906\,631) & = \phantom{-}3579, \displaybreak[0] \\
R (81\,916\,669\,933\,751) & = -2397,           \displaybreak[0] \\
R (83\,575\,544\,636\,251) & = -1884,           \displaybreak[0] \\
R (84\,688\,857\,018\,011) & = -1183,           \displaybreak[0] \\
R (84\,771\,692\,838\,421) & = \phantom{-}2281,                  \\
R (85\,325\,902\,236\,661) & = -4473.                            
\end{align*}
There have been
$5\,235$
positive and
$5\,557$
negative quotients detected.
\end{ttt}

\begin{rems}
\begin{abc}
\item
We note explicitly that this is not at all a constructed example. One may basically consider every interval which is not too small and will observe the same phenomena.
\item
Being very sceptical one might raise the objection that the computations done in our program do not really prove that the
$10\,792$
numbers~$p$
which appear in the data are indeed prime.

It is, however, very unlikely that one of them is composite as they all passed two~tests. First, they passed the sieve which in this case makes sure they have no prime
divisor~$\leq 8\,302\,871$.
This means, if one is composite then it decomposes into the product of two almost equally large primes. Furthermore, they were all found probably~prime by the Fibonacci composedness
test~$p | F_{p-1}$.

It is easy to check primality for all of them by a separate program.
\end{abc}
\end{rems}

\begin{ttt}
{\bf Statistical Observations.}
A more spectacular interval
is~$[0, 10^{12}]$.
One may expect a lot more small quotients as all small prime numbers are taken into~consideration.

Here, we may do some statistical analysis on the small positive values
of the quotient~$R^\prime$
which is given by
$R^\prime (p) := (F_{p - 1} \mod p^2) / p$
for
$p \equiv \pm1 \pmod 5$
and by
$R^\prime (p) := (F_{2p + 2} \mod p^2) / p$
for
$p \equiv \pm2 \pmod 5$.

Our computations show that there exist
$96\,909$~quotients
less~than~$100\,000$,
$12\,162$~quotients
less~than~$10\,000$,
$1\,580$~quotients
less~than~$1\,000$,
$216$~quotients
less~than~$100$,
and
$30$~quotients
less~than~$10$.
The latter ones are

\begin{quell}
3 ist eine interessante Primzahl. Quotient 1
7 ist eine interessante Primzahl. Quotient 1
11 ist eine interessante Primzahl. Quotient 5
13 ist eine interessante Primzahl. Quotient 7
17 ist eine interessante Primzahl. Quotient 2
19 ist eine interessante Primzahl. Quotient 3
43 ist eine interessante Primzahl. Quotient 8
89 ist eine interessante Primzahl. Quotient 5
163 ist eine interessante Primzahl. Quotient 6
199 ist eine interessante Primzahl. Quotient 5
239 ist eine interessante Primzahl. Quotient 5
701 ist eine interessante Primzahl. Quotient 5
941 ist eine interessante Primzahl. Quotient 6
997 ist eine interessante Primzahl. Quotient 3
1063 ist eine interessante Primzahl. Quotient 2
1621 ist eine interessante Primzahl. Quotient 2
2003 ist eine interessante Primzahl. Quotient 1
27191 ist eine interessante Primzahl. Quotient 8
86813 ist eine interessante Primzahl. Quotient 6
123863 ist eine interessante Primzahl. Quotient 2
199457 ist eine interessante Primzahl. Quotient 7
508771 ist eine interessante Primzahl. Quotient 2
956569 ist eine interessante Primzahl. Quotient 4
1395263 ist eine interessante Primzahl. Quotient 3
1677209 ist eine interessante Primzahl. Quotient 1
3194629 ist eine interessante Primzahl. Quotient 5
11634179 ist eine interessante Primzahl. Quotient 2
467335159 ist eine interessante Primzahl. Quotient 4
1041968177 ist eine interessante Primzahl. Quotient 6
6_71661_90593 ist eine interessante Primzahl. Quotient 1
\end{quell}

Except for
$0$
and
$9$,
all one-digit numbers do appear.

Further, the counts are again well within the expected range. For example, consider one-digit numbers.
$R(3)$
and
$R(7)$
are automatically one-digit. Therefore, the expected count is
$$2 + \sum_{\hsmash{\substack{p = 10 \\ p \rmtext{ prime}}}}^{10^{12}} \frac{10}p \approx 2 + 10 \cdot (\log (\log 10^{12}) - \log (\log 10)) \approx 26.849\,066$$
which is surprisingly close
the~$30$
one-digit quotients which were actually~found.

We note that already for two-digit quotients, it is no~longer true that they appear only within the
subinterval~$[0, 10^{11}]$.
In fact, there are twelve prime
numbers~$p \in [10^{11}, 10^{12}]$
such that
$R^\prime (p) < 100$.
These are the following.\smallskip\pagebreak[3]

\begin{quell}
101876918491 liefert 87
115301883659 liefert 60
129316722167 liefert 44
147486235177 liefert 59
170273590301 liefert 78
233642484991 liefert 89
261836442223 liefert 45
277764184829 liefert 64
283750593739 liefert 37
305128713503 liefert 93
334015396151 liefert 79
442650398821 liefert 74
\end{quell}

Once again, we may compare this to the expected count which is here
$$\sum_{\hsmash{\substack{p = 10^{11} \\ p \rmtext{ prime}}}}^{10^{12}} \frac{100}p \approx 100 \cdot (\log(\log 10^{12}) - \log(\log 10^{11})) \approx 8.701\,137\,73.$$
\end{ttt}

\renewcommand{\thefootnote}{\fnsymbol{footnote}}
\footnotetext[0]{version of Decembre 1, 2004}

\frenchspacing
\begin{thebibliography}{Wa}
\bibitem[Co]{Cohen}
Cohen, H.: A course in computational algebraic number theory, {\em Sprin\-ger,} Graduate Texts Math.~138, Berlin~1993
\bibitem[Fo]{Forster}
Forster, O.: Algorithmische Zahlentheorie (Algorithmic number theory), {\em Vie\-weg}, Braunschweig~1996
\bibitem[G\"o]{Gottsch}
G\"ottsch, G.: \"Uber die mittlere Periodenl\"ange der Fibonacci-Folgen modulo~$p$\; (On~the average period length of the Fibonacci sequences modulo~$p$), {\em Disser\-tation,} Fakult\"at f\"ur Math.\ und Nat.-Wiss., Hannover~1982
\bibitem[IR]{IR}
Ireland, K., Rosen, M.: A classical introduction to modern number theory,
Second edition, {\em Sprin\-ger,} Graduate Texts Math.~84, New~York~1990
\bibitem[Ja]{Jarden}	
Jarden, D.: Two theorems on Fibonacci's sequence, {\em Amer.\ Math.\ Monthly\/} 53\br(1946)\br425--427
\bibitem[Mo]{Montgomery}
Montgomery, P.~L.: Modular multiplication without trial division, {\em Math.\ Comp.} 44\br(1985)\br519--521
\bibitem[RS]{RS}
Rosser, J.~B., Schoenfeld, L.: Approximate formulas for some functions of prime numbers, {\em Illinois J.\ Math.} 6\br(1962)\br64--94
\bibitem[Wa]{Wall}
Wall, D.~D.: Fibonacci series modulo~$m$, {\em Amer.\ Math.\ Monthly\/} 67\br(1960)\br525--532
\end{thebibliography}
\end{document}